\documentclass[12pt]{article}
\usepackage{amsmath}[1996/11/01]
\usepackage{amssymb,amsthm,amsfonts,latexsym,psfig,epsfig,graphics}
\usepackage{amscd}
\usepackage{epic}
\advance \textheight by -1 \topmargin
\advance \textheight by -1 \headheight
\advance \textheight by -1 \headsep
\advance \textheight by -2 \footskip
\vsize = \textheight
\hsize = \textwidth
\newtheorem{theorem}{Theorem}[section]
\newtheorem{satz}[theorem]{Theorem}
\newtheorem{notation}[theorem]{Notation}
\newtheorem{prop}[theorem]{Proposition}

\newtheorem{kor}[theorem]{Corollary}

\newcommand{\bew}{\noindent\underline{Proof.}\ }

\newtheorem{rem}[theorem]{Remark}
\newtheorem{lemma}[theorem]{Lemma}

\newtheorem{definition}[theorem]{Definition}
\newtheorem{folg}[theorem]{Corollary}

\newcommand{\fa}{\mbox{ for}\mbox{ all }}
\newcommand{\be}{\begin{enumerate}}
\newcommand{\ee}{\end{enumerate}}
\newcommand{\bi}{\begin{itemize}}
\newcommand{\ei}{\end{itemize}}

\newcommand{\ba}{\begin{array}}
\newcommand{\ea}{\end{array}}

\newcommand{\ra}{\rightarrow}

\newcommand{\tr}{\mbox{trace}}

\newcommand{\Z}{{\mathbb{Z}}}
\newcommand{\Q}{{\mathbb{Q}}}

\newcommand{\N}{{\mathbb{N}}}

\newcommand{\eb}{\phantom{zzz}\hfill{$\square $}\smallskip}

\renewcommand{\em}{\sf}

\DeclareMathOperator{\Hom}{Hom}
\DeclareMathOperator{\rad}{rad}

\DeclareMathOperator{\Id}{Id}
\DeclareMathOperator{\Tr}{Tr}

\begin{document}

\Huge
\begin{center}
{\bf On the radical idealizer chain of symmetric orders.}
\end{center}
\normalsize
\begin{center}
Gabriele Nebe \footnote{Lehrstuhl D f\"ur Mathematik, RWTH Aachen,
52056 Aachen, nebe@math.rwth-aachen.de}
\end{center}

\small

{\sc Abstract}:
If $\Lambda $ is an indecomposable, non-maximal,
 symmetric order, then the idealizer of the
radical $\Gamma := \Id(J(\Lambda )) = J(\Lambda )^{\#} $ is the dual of the
radical. If $\Gamma $ is  hereditary 
then $\Lambda $ has a Brauer tree
(under modest additional assumptions).
Otherwise $\Delta := \Id(J(\Gamma )) = (J(\Gamma )^2)^{\#} $.
If $\Lambda = \Z_p G$ for a $p$-group $G\neq 1$, then $\Gamma $ is hereditary 
iff $G\cong C_p$ and otherwise  $[\Delta : \Lambda ] = p^2 | G/(G'G^p)| $.
For Abelian groups $G$, the length of the radical idealizer chain 
of $\Z_pG$ is $(n-a)(p^{a} - p^{a-1})+p^{a-1}$, where $p^n$ is the order and
$p^a$ the exponent of the Sylow $p$-subgroup of $G$.

\normalsize 

\section{Introduction}
Throughout this paper 
let $R$ be a discrete valuation ring with maximal ideal 
$\pi R$  and residue class field $R/\pi R =: k$.
Let $K$ be the field of fractions of $R$ and 
${\cal A}$  a separable finite dimensional $K$-algebra.
If $\Lambda $ is an $R$-order in ${\cal A}$, then there is a canonical
process, the so called {\em radical idealizer process},
that constructs an ascending chain of over-orders of $\Lambda $ ending in 
a hereditary order $\Lambda _N $, called the {\em head order}
of $\Lambda $ (see Remark \ref{radid}).
We call the length of this chain the {\em radical idealizer length}
$l_{rad}(\Lambda )$ of $\Lambda $.
Hereditary orders are well understood, see \cite{Jac}, \cite[Chapter 9]{Rei}.
 They are direct sums of hereditary orders in the simple components 
of ${\cal A}$.
So one might hope to classify $R$-orders according to the length 
of the radical idealizer chain and the head order.

An important class of $R$-orders are the {\em symmetric orders},
that are self dual with respect to some trace bilinear form on ${\cal A}$.
Examples of symmetric orders are provided by 
blocks of group rings  $RG$ for finite groups $G$.
The main tool to deal with symmetric orders is Jacobinski's conductor formula
(see Theorem \ref{fuehrer}) stating that for any  over-order $\Gamma $ of a symmetric
order $\Lambda $ the conductor $F_{\Gamma }(\Lambda )$ (which is the largest 
$\Gamma $ ideal in $\Lambda $) is the dual of $\Gamma $.
If $\Gamma  = \Id (J(\Lambda ))$ is the idealizer of $J(\Lambda )$, then 
a converse of this formula holds:
Theorem  \ref{jistdual} shows that 
for indecomposable, non-hereditary, symmetric orders $\Lambda $ the 
dual of $J(\Lambda )$ is the idealizer of $J(\Lambda )$.
Using his conductor formula, Jacobinski shows that the indecomposable
symmetric orders $\Lambda $ with $l_{rad}(\Lambda ) = 0$ 
are maximal orders (Theorem \ref{erbsymmax}).
If $l_{\rad }(\Lambda ) = 1$, then 
 one may derive the Brauer tree of $\Lambda $ using the
idealizer of $J(\Lambda )$ (see Proposition \ref{tree} and 
 \cite[Section 11]{Jac}).
In the present paper the first two steps of the radical idealizer chain for
symmetric orders are investigated and properties of  symmetric orders 
$\Lambda $ with $l_{rad}(\Lambda ) = 2$ are determined.
We apply the theorems to $p$-groups $G$ showing that 
$l_{rad}(\Z_pG) = 1$ if and only if $G \cong C_p$,
$l_{rad}(\Z_pG ) = 2 $ if and only if $G \cong C_4$ or $G\cong C_2\times C_2$.
Moreover we calculate $l_{rad}(\Z_pG)$ for Abelian groups $G$.

\section{The radical idealizer chain}

Let $\Lambda $ be an $R$-order in ${\cal A}$.
Then the {\em Jacobson radical} $J(\Lambda )$ is 
the intersection of all maximal right ideals of $\Lambda $.
It is a 2-sided ideal of $\Lambda $, in fact the 
smallest ideal $I$ of $\Lambda $, such that $\Lambda /I$ is
a semi-simple $k$-algebra.
One other important characterization of $J(\Lambda )$ is
that $J(\Lambda )$ is the biggest $\Lambda $-ideal $I$ in $\Lambda $ that
is {\em pro-nilpotent}, i.e. for which there is $m\in \N$ such 
that $I^m \subset \pi \Lambda $ (cf. \cite[Lemma 8.5]{Jac}).

\begin{definition} (see \cite[Section 39]{Rei})
Let $\Lambda $, $\Lambda '$ be $R$-orders in ${\cal A}$.
Then $\Lambda $ {\em radically covers} $\Lambda '$, $\Lambda \succ \Lambda '$,
if $\Lambda \supseteq \Lambda '$ and $J(\Lambda ) \supseteq J(\Lambda ')$.
If $\Lambda $ is maximal with respect to $\succ $, then $\Lambda $ is
called {\em extremal}.
\end{definition}

\begin{lemma}
Let $\Gamma \succ \Lambda $ be two $R$-orders in ${\cal A}$.
Then $J(\Lambda ) = J(\Gamma ) \cap \Lambda $ and 
$\Lambda /J(\Lambda ) $ is isomorphic to a sub-algebra of 
$\Gamma /J(\Gamma )$.
Moreover every simple $\Gamma $-module is semi-simple as a $\Lambda $-module.
\end{lemma}

\bew
Since $J(\Gamma ) \cap \Lambda $ is an  ideal of $\Lambda $ that is 
nilpotent modulo $\pi \Lambda $,
it is contained in $J(\Lambda )$.
On the other hand $J(\Lambda ) \subseteq J(\Gamma )$, because 
$\Gamma \succ \Lambda $.
Therefore $J(\Gamma )\cap \Lambda = J(\Lambda )$ and 
$\Lambda /J(\Lambda ) \cong (\Lambda +J(\Gamma ))/J(\Gamma )$ 
is naturally embedded in $\Gamma /J(\Gamma )$.
The second assertion follows from the fact that
 $\Gamma J(\Lambda )=  J(\Lambda ) \subseteq J(\Gamma )$
 which implies that $\Gamma /J(\Gamma )$ is a semi-simple $\Lambda $-module.
\eb


Recall that an order $\Gamma $ is called {\em hereditary}, if 
every left ideal of $\Gamma $ is projective (see \cite[Section 10]{Rei}). 

\begin{satz}(\cite[Satz 8.12]{Jac}, \cite[Theorem 39.14]{Rei})
An $R$-order $\Lambda $ in ${\cal A}$ is extremal, if and only if
$\Lambda $ is hereditary.
\end{satz}

\begin{definition}
Let $L$ be a full $R$-lattice in ${\cal A}$.
The {\em left order} of $L$ is $O_l(L):=\{ a\in {\cal A} \mid aL\subseteq L \}$.
Analogously one defines the {\em right order} $O_r(L)$ of $L$.
$\Id(L) := O_l(L)\cap O_r(L)$ is called the {\em idealizer} of $L$.
\end{definition}

\begin{rem}
Let $\Lambda $ be an order and let  $\Gamma $ be one of
$O_l(J(\Lambda ))$, $O_r(J(\Lambda )$, or $\Id(J(\Lambda ))$. 
Then $\Gamma \succ \Lambda $.
\end{rem}

The following characterization of hereditary orders 
 is shown in \cite[Theorem 39.11]{Rei} for $O_l(J(\Lambda ))$ instead of $\Id
(J(\Lambda) )$. 
With a completely analogous proof (see \cite{habil}) one shows

\begin{satz}
Let $\Lambda $ be an $R$-order in ${\cal A}$.
Then $\Lambda  = \Id(J(\Lambda )) $ if and only if $\Lambda $ is hereditary.
\end{satz}

\begin{rem}(cf. \cite{BeZ}){\label{radid}}
Letting $\Lambda _0:=\Lambda $ and $\Lambda _{n+1}:=\Id(J(\Lambda _n))$
for $n=0,1,2,\ldots $ defines a canonical process,
the so called {\em radical idealizer process} that constructs from 
an $R$-order $\Lambda $ in ${\cal A}$ successively bigger $R$-orders 
$\Lambda _0 \subset \Lambda _1 \subset \ldots \subset \Lambda _N = \Lambda
_{N+1}$, the so called {\em radical idealizer chain}. 
The order $\Lambda _N$ is hereditary and called the
{\em head order} of $\Lambda $.
If $N$ is minimal such that $\Lambda _N = \Lambda _{N+1}$, then 
$N$ is called the {\em radical idealizer length} $l_{rad}(\Lambda )$
of $\Lambda $.
\end{rem}

Replacing $\Id $ by $O_l$ respectively $O_r$, one can define
left- and right-idealizer chain similarly.
Since automorphisms of $\Lambda $ preserve the radical,
they also yield automorphisms of $O_r(J(\Lambda ))$, $O_l(J(\Lambda ))$
and $\Id (J(\Lambda ))$.
The advantage of taking two-sided idealizers is that
$\Id (J(\Lambda ))$ is also preserved under anti-automorphisms of
$\Lambda $, which interchange $O_r(J(\Lambda ))$ and $O_l(J(\Lambda ))$
and hence left- and right-idealizer chains.

The next remark gives a lower bound on the length of the radical idealizer chain
and is also useful for the explicit calculation of $\Id(J(\Lambda ))$,
since one may calculate modulo the maximal ideal $\pi R$:

\begin{rem}{\label{piinJ}}
Let $\Lambda $ be an $R$-order in ${\cal A} $ and 
$\Gamma := \Id(J(\Lambda ))$. Then 
$J(Z(\Lambda )) \Gamma \subseteq \Lambda $, in particular
$\pi \Gamma \subseteq \Lambda $.
\end{rem}

\bew
$J(Z(\Lambda )) \Lambda \subset J(\Lambda ) \subset \Lambda $, 
since $J(Z(\Lambda ))\Lambda $ is
nilpotent modulo  $\pi \Lambda $.
Therefore $$J(Z(\Lambda )) \Gamma = J(Z(\Lambda )) \Lambda \Gamma \subseteq J(\Lambda ) \Gamma =
J(\Lambda ) \subseteq \Lambda .$$
\eb

\begin{definition}
Let $\Lambda $ be an $R$-order in ${\cal A}$ and let 
$\epsilon _1,\ldots ,\epsilon _s$ be the central primitive idempotents
of ${\cal A}$. 
Then the {\em defect} of $\Lambda $ is the minimal $d$ 
such that $\pi ^d \epsilon _t \in \Lambda $ for all $1\leq t \leq s$.
\end{definition}

Note that this coincides with the usual definition of defect
for blocks of group rings,
if $K$ is an unramified extension of $\Q _p$.

Since hereditary orders contain the central primitive idempotents of 
${\cal A}$, one gets the following corollary:

\begin{kor}
The radical idealizer length $l_{rad}(\Lambda )$ is greater or equal
than the defect of $\Lambda $.
\end{kor}

\begin{definition}
For two $R$-orders $\Lambda $, $\Gamma $ in ${\cal A}$, the 
{\em conductor} of $\Gamma $ in $\Lambda $ is the biggest 
$2$-sided $\Gamma $-ideal $F_{\Gamma }(\Lambda )$ that is contained in $\Lambda $.
Analogously one defines 
the {\em left conductor} $F^{(l)}_{\Gamma }(\Lambda )$ and 
the {\em right conductor} $F^{(r)}_{\Gamma }(\Lambda )$ as the
largest left- respectively right-ideal of $\Gamma $ contained in $\Lambda $.
\end{definition}

The following lemma is a straightforward generalization of 
\cite[Theorem 2.2]{CPW}:

\begin{lemma}\label{commut}
Assume that  ${\cal A}$ is commutative,
let $\epsilon _1,\ldots , \epsilon _s$
be the primitive idempotents in ${\cal A}$ and assume that 
$\Gamma := \oplus _{i=1}^s \epsilon _i \Lambda $ is the maximal order in ${\cal A}$.
For $i\in \{ 1,\ldots , s\}$ let $\pi _i $ be a prime element in $\epsilon _i\Lambda $and put $\pi := (\pi _1 ,\ldots , \pi _s) \in \Gamma $.
Let $$\Lambda = \Lambda _0 \subset \Lambda _1 \subset \ldots \subset \Lambda _N  = \Gamma $$
be the radical idealizer chain of $\Lambda $. 
Then for $n=0,\ldots , N$ 
$$F_{\Gamma }(\Lambda _n ) = \pi ^{-n} F_{\Gamma }(\Lambda ) \cap \Gamma . $$
\end{lemma}

\proof
We argue by induction on $n$, the case $n=0$ being trivial.
Assume that 
$$F_{\Gamma }(\Lambda _n ) = \pi ^{-n} F_{\Gamma }(\Lambda ) \cap \Gamma 
= \bigoplus _{i=1}^s \pi _i ^{a_i} \Lambda \epsilon _i . $$
Splitting off the direct summands of $\Lambda _n$ that are maximal orders,
we may assume that $a_i > 0$ for all $i$. 
Then $$ F_{\Gamma } (\Lambda _n) \subseteq J(\Lambda _n) = \pi \Gamma \cap \Lambda _n .$$
Since $F_{\Gamma }(\Lambda _n)$ is a $\Gamma $-ideal one gets 
$$\pi ^{-1} F_{\Gamma }(\Lambda _n) J(\Lambda _n) \subseteq 
\pi ^{-1} F_{\Gamma }(\Lambda _n) \pi \Gamma \subseteq 
F_{\Gamma }(\Lambda _n) \subseteq  J(\Lambda _n ) .$$
Hence $$\pi ^{-1} F_{\Gamma }(\Lambda _n) \subseteq \Id (J(\Lambda _n)) = \Lambda _{n+1}$$ and therefore 
$$F_{\Gamma }(\Lambda _{n+1} ) \supseteq \pi ^{-(n+1)} F_{\Gamma }(\Lambda ) \cap \Gamma  .$$
The opposite inclusion follows from Remark \ref{piinJ}.
\eb

\begin{folg}\label{folgcommut}
In the notation of Lemma \ref{commut} let
$$F_{\Gamma }(\Lambda ) = \bigoplus _{i=1}^s \pi_i ^{a_i}  \Lambda \epsilon _i .$$
Then $l_{rad}(\Lambda ) = \max _{i=1,\ldots , s} a_i $.
\end{folg}

\section{Idealizers and bilinear forms}

Idealizers can be calculated using a symmetric non-degenerate
bilinear form  $$\phi : {\cal A} \times {\cal A} \ra K \mbox{ that is 
{\em associative}, i.e. } \phi (ab,c) = \phi (a,bc) \mbox{ for all } a,b,c \in 
{\cal A}.$$
It is easy to see that such an associative
bilinear form $\phi $ is of the form
$$\Tr_z: {\cal A} \times {\cal A} \ra K , (a,b) \mapsto tr_{red} (zab) ,$$
where $z\in Z({\cal A})^*$ is an invertible element of the center of ${\cal A}$
and $tr_{red}$ denotes the reduced trace of ${\cal A}$.
Fix such an associative symmetric bilinear form $\phi = \Tr_z$.
For a full $R$-lattice $L$ in ${\cal A}$ let 
$$L^{\#} := \{ a\in {\cal A} \mid \phi(L,a)\subset R \} $$
be the {\em dual lattice} with respect to $\phi $.
It is frequently used that dualizing is an inclusion reversing
bijection of the set of full $R$-lattices in ${\cal A}$ and
that $(L^{\#})^{\#} = L$ for all full $R$-lattices $L$ in ${\cal A}$.

\begin{lemma}{\label{dualistid}}
Let $\Gamma $ be an $R$-order in ${\cal A}$. 
Then $\Gamma ^{\#}$ is a $2$-sided $\Gamma $-ideal with 
$$\Gamma = O_l(\Gamma ^{\#} ) = O_r(\Gamma ^{\#})  =\Id (\Gamma ^{\#}) .$$
\end{lemma}

\bew
Let $x,\gamma \in \Gamma $ and $y\in \Gamma ^{\# }$.
Then $\phi(x,\gamma y) = \phi (x\gamma,y) \in R$ and
$\phi(y\gamma,x) = \phi(y,\gamma x) \in R$ and therefore $\Gamma \subset \Id(\Gamma ^{\#})$.
On the other hand let $\lambda \in O_l(\Gamma ^{\# })$.
Then for all $y\in \Gamma ^{\#} $ and 
 $x\in \Gamma = (\Gamma ^{\#} )^{\#} $
$$\phi (x\lambda , y ) = \phi(x,\lambda y ) \in R .$$
Hence $\Gamma \lambda \subseteq (\Gamma ^{\#})^{\#} = \Gamma $ and 
therefore $\lambda = 1 \lambda \in \Gamma $.
Analogously one gets $O_r(\Gamma ^{\# }) \subseteq \Gamma $.
\eb

\begin{prop}{\label{ddd}}
If $L$ is a full $R$-lattice in ${\cal A}$, then 
$O_l(L) = (LL^{\#})^{\#} $, $O_r(L)=(L^{\#}L)^{\#}$, and hence
$$\Id(L) = (LL^{\#})^{\#} \cap (L^{\#}L)^{\#} .$$
\end{prop}

\bew
We only show $O_l(L) = (LL^{\#})^{\#} $. 
Let $\gamma \in {\cal A}$.
Then $\gamma \in (LL^{\# })^{\#} $ if and only if
 for all $x\in L$, $y\in L^{\#}$ 
$$\phi (\gamma ,xy) = \phi (\gamma x,y) \in R$$ 
which is equivalent to $\gamma L \subseteq (L^{\#})^{\#}  = L$, i.e.
$\gamma \in O_l(L)$.
Analogously 
$O_r(L) = (L^{\#}L)^{\#}$.
\eb

From this one gets an interesting direct description of 
the idealizer of the radical.

\begin{kor} 
Let $\Lambda $ be an $R$-order in ${\cal A}$ and 
$\Gamma := \Id(J(\Lambda ))$.
Then $\Gamma $ is the biggest $\Lambda $-ideal $I \subset \frac{1}{\pi }\Lambda $
such that $I/J(\Lambda )$ is a semi-simple $\Lambda$-$\Lambda $-bimodule.
\end{kor}

\bew
By Proposition \ref{ddd}
$$\Gamma ^{\#} = (J(\Lambda ) J(\Lambda )^{\#}) + 
(J(\Lambda )^{\#} J(\Lambda )) $$
is the smallest $\Lambda $-ideal $J$ in $J(\Lambda )^{\# }$ for which
$J(\Lambda )^{\#} / J$ is a semi-simple $\Lambda $-$\Lambda $-bimodule.
Since the dual of a bimodule is semi-simple if and only if the module
is semi-simple, the corollary follows.
\eb

\section{Symmetric orders.}

\begin{definition}
An $R$-order $\Lambda $ in ${\cal A}$ is called {\em symmetric},
if there is a non-degenerate, symmetric, associative, bilinear form
$\phi = \Tr _z : {\cal A} \times {\cal A} \ra K$ with 
$\Lambda = \Lambda ^{\# }$.
\end{definition}

\begin{lemma}(see e.g. \cite[Proposition (1.6.2)]{The}){\label{efd=fe}}
Let $e,f \in \Lambda $ be two idempotents in the symmetric order $\Lambda $.
Then $\phi _{|e\Lambda f \times f \Lambda e}$ is a regular
$R$-bilinear pairing.
In particular $e\Lambda e$ is a symmetric order in $e{\cal A} e$.
\end{lemma}

An important tool to deal with symmetric orders is Jacobinski's 
conductor formula:

\begin{satz}(\cite[Satz 10.6]{Jac}){\label{fuehrer}}
Let $\Lambda $ be a 
 symmetric $R$-order in ${\cal A}$ and 
$\Gamma \supseteq \Lambda $ an over-order. 
Then left and right conductor coincide and are equal to the dual of $\Gamma $:
$$F_{\Gamma }(\Lambda ) =F_{\Gamma }^{(l)}(\Lambda )
 =F_{\Gamma }^{(r)}(\Lambda ) = \Gamma ^{\# }. $$
\end{satz}

\bew
Since 
$\Gamma ^{\# } \subseteq \Lambda $ is a $\Gamma $-ideal by 
Lemma \ref{dualistid}
one has $ \Gamma ^{\# } \subseteq F_{\Gamma }(\Lambda) $.
On  the other hand  if
 $x\in F_{\Gamma }(\Lambda )$ and $\gamma \in \Gamma $  then
$\phi (x,\gamma) = \phi(x\gamma , 1) \in R$ is integral, since 
$x\gamma \in \Lambda = \Lambda ^{\#}$.
\eb

From this proof one even gets that 
$\Gamma ^{\# }$  is the largest $R$-lattice $L$ in $\Lambda $ with
$\Gamma L \subseteq \Lambda $.

\cite[Satz 10.7]{Jac} and \cite[Theorem III.8]{Ple 83}
describe the conductor of $F_{\Gamma }(\Lambda ) $ for hereditary 
and (more general) graduated over-orders $\Gamma $ of  the symmetric 
order $\Lambda $.
To apply this precise version of the conductor formula, we need the following
(technical) notation:

\begin{notation}\label{condprec}
Let $\epsilon _1,\ldots , \epsilon _s$ be the central primitive idempotents
in ${\cal A}$.
\\
Let
$z\in Z({\cal A})$ be such that
$\Lambda = \Lambda ^{\# }$ with respect to
$Tr_z$ and
\\
$z_i=\epsilon_i z \in Z(\epsilon _i {\cal A}) = K _i$ ($1\leq i\leq s$).
\\
Let $R_i$ be the maximal order in $K_i$ with maximal ideal
$\wp _i$,
\\
$\wp _i^{-d_i}$ the inverse different of $R_i$ over $R$
and $n_i \in \Z $ with $z_i R_i = \wp_i^{-n_i} $
($1\leq i \leq s$).
\\
The simple algebra ${\cal A} \epsilon _i$ is isomorphic to a matrix
ring over a central $K_i$-division algebra $D_i$.
Let $\Omega _i$ be the maximal order in $D_i$ and $m_i^2  = \dim _{K_i}(D_i)$.
\end{notation}

\begin{theorem}( \cite[Satz 10.7]{Jac}, \cite[Theorem III.8]{Ple 83}) \label{condform}
With the notation above let $\Delta $ be a hereditary order in ${\cal A}$.
Then 
$$\Delta ^{\#}  = \oplus _{i=1}^s \wp _i ^{m_i(n_i-d_i-1)} J(\Delta ) \epsilon _i .$$
\end{theorem}

\begin{satz}{\label{erbsymmax}}
Let $\Lambda $ be an indecomposable  symmetric $R$-order in ${\cal A}$
and  $0\neq e^2=e \in \Lambda $ be an idempotent such that $e\Lambda e$
is hereditary. 
Then $\Lambda $ is a maximal order.
\end{satz}

\bew
$\epsilon _i e \Lambda e =: \Lambda _i$ is either $\{ 0 \}$ 
or a symmetric hereditary order in $e {\cal A } e $ for all $1\leq i \leq s $.

Let $i$ be fixed such that $\Lambda _i \neq \{ 0 \}$.
Then the conductor formula \ref{condform}
yields that 
$$\Lambda _i = \Lambda_i^{\#} = \wp _i ^{m_i(n_i-d_i-1)}  J(\Lambda _i) .$$
In particular $\Lambda _i$ is isomorphic to $J(\Lambda _i)$ as a bimodule 
and therefore $\Lambda_i$ a maximal order in $\epsilon _i e {\cal A} e$,
$m_i = 1$  and
$n_i = d_i$.
But then the conductor of every maximal over-order $\Gamma $ of $\Lambda $ 
in $\Lambda $ is of the form $F_{\Gamma }(\Lambda ) = \Gamma ' \oplus 
\epsilon _i \Gamma $ for a suitable order $\Gamma '$.
In particular $\epsilon _i \in \Lambda $.
Since $\Lambda $ is indecomposable $\Lambda = \epsilon _i \Lambda = \epsilon _i
\Gamma $ and $\Lambda $ is a maximal order in the simple 
$K$-algebra ${\cal A} = \epsilon _i {\cal A}$. 
\eb

Putting $e=1$ in
Theorem \ref{erbsymmax} this characterizes the symmetric orders 
$\Lambda $ with $l_{rad} (\Lambda ) = 0$ as
maximal orders. In particular if 
$\Lambda $ is a block of a group ring $RG$, then 
$l_{rad} (\Lambda ) = 0 $ if and only if the defect of
$\Lambda $ is 0 (see \cite[Satz 11.1]{Jac}).

\section{The radical idealizer of symmetric orders.}

In this and the next section the first two steps of the radical idealizer chain
of symmetric orders are made precise.
The first theorem is a sort of converse of the conductor formula.

\begin{satz}{\label{jistdual}}
Let $\Lambda $ be a non-hereditary, indecomposable, symmetric $R$-order in
${\cal A}$ and $\Gamma := \Id(J(\Lambda ))$ the idealizer of the radical of 
$\Lambda $. 
Then $\Gamma = J(\Lambda )^{\# }$.
\end{satz}

\bew
$\Gamma ^{\#} \subseteq \Lambda $ is the largest $\Gamma $-ideal in $\Lambda $
by the conductor formula. 
Since $J(\Lambda ) \subseteq \Lambda $ is a $\Gamma $-ideal, one has
$J(\Lambda ) \subseteq \Gamma ^{\#} $ and therefore $\Gamma \subseteq J(\Lambda )^{\# }$.

To show the converse inclusion let
$e_1,\ldots , e_h \in \Lambda $ be orthogonal 
idempotents that map onto the central primitive
idempotents of $\Lambda /J(\Lambda )$  with $1=e_1+\ldots + e_h $.
Then $\Lambda = \oplus _{i,j=1}^h e_i\Lambda e_j$, 
$(e_i\Lambda e_j)^{\# } = e_j \Lambda e_i$ and
$e_i\Lambda e_i$ is a symmetric $R$-order in $e_i {\cal A} e_i$,
which is not hereditary because of Theorem 
 \ref{erbsymmax}.
Now
$$J(\Lambda ) = \bigoplus _{i\neq j = 1}^h e_i \Lambda e_j \oplus
\bigoplus _{i=1}^h J(e_i\Lambda e_i)$$
and with Lemma \ref{efd=fe}
$$J(\Lambda )^{\# } = \bigoplus _{i\neq j = 1}^h e_i \Lambda e_j \oplus
\bigoplus _{i=1}^h J(e_i\Lambda e_i)^{\# }.$$

Assume first that $h=1$ so $\Lambda /J(\Lambda )$ is a simple $k$-algebra.
Then $J(\Lambda )$ is a maximal $2$-sided ideal in $\Lambda $.
Therefore either $\Gamma ^{\# } = \Lambda $ or $\Gamma ^{\#} = J(\Lambda )$.
In the first case $\Gamma = (\Gamma ^{\# })^{\#}  = \Lambda ^{\# }
= \Lambda $ and therefore $\Lambda $ is hereditary contradicting the
assumption. In the second case $\Gamma = J(\Lambda )^{\# }$ 
and the theorem follows.

Now let $h$ be arbitrary and $\Gamma _j := J(e_j\Lambda e_j)^{\# }$.
From above $\Gamma _j = \Id(J(e_j\Lambda e_j))$  is an order, so 
it remains to show that for
 $i\neq j$ the summand
$e_i\Lambda e_j$ of $J(\Lambda )^{\#}$ is a  $\Gamma _i$-$\Gamma _j$-bimodule.
The inclusion $(e_i\Lambda e_j)\Gamma _j \subseteq e_i \Lambda e_j$
is equivalent to 
$$e_j \Lambda e_i = (e_i\Lambda e_j)^{\#} \subseteq
(e_i \Lambda e_j \Gamma _j)^{\#}.$$
So let $e_j a e_j \in \Gamma _j$ with  $a\in J(\Lambda )^{\#}$ and
$\gamma , \lambda \in \Lambda $.
Then
$$\phi (e_i \lambda e_j e_j a e_j ,e_j \gamma e_i ) =
\phi (e_j \gamma e_i , e_i \lambda e_j a e_j)
=\phi (e_j \gamma e_i  \lambda e_j , e_j a e_j)
\in R,$$ because
$e_j\gamma e_i \lambda e_j \in J(\Lambda )$ (note that $i\neq j$)
and $a\in J(\Lambda )^{\#}$.
Analogously $\Gamma _i (e_i \Lambda e_j)\subset e_i \Lambda e_j$.

Since $J(e_j\Lambda e_j)$ is a  $\Gamma _j$-bimodule ($1\leq j \leq s$),
one gets  $J(\Lambda )^{\#} J(\Lambda ) J(\Lambda )^{\# } \subseteq
J(\Lambda )$ and hence $J(\Lambda )^{\# } \subseteq \Gamma $.
\eb

Since $\Id (J(\Lambda ) ) \subseteq O_r(J(\Lambda ))$ and
the right-conductor $O_r(J(\Lambda ))^{\#} \supset J(\Lambda ) $
one gets the same result for the left- and right-idealizer of
$J(\Lambda )$.

\begin{kor}{\label{l=r}}
Let $\Lambda$ be a symmetric order.
Then $$\Id (J(\Lambda )) = O_r (J(\Lambda )) = O_l (J(\Lambda )) .$$
\end{kor}

The orders $e \Lambda e$, where $e^2=e \in \Lambda $ is an idempotent
in $\Lambda $ mapping onto a central primitive idempotent of 
$\Lambda /J(\Lambda )$ play an important role in the above proof.
These orders are {\em 2-sided local} orders, i.e. they have a
unique maximal 2-sided ideal. These orders have a unique simple module,
or equivalently $e\Lambda e/eJ(\Lambda ) e$ is a simple $k$-algebra.

\begin{lemma}{\label{selbstduallok}}
Let $\Lambda $ be a 2-sided local, symmetric $R$-order and
$\Gamma := \Id(J(\Lambda ))$.
Then either $\Gamma /J(\Gamma ) \cong \Lambda /J(\Lambda )$
as $\Lambda $-$\Lambda $-bimodule and $J(\Gamma )= J(\Gamma )^{\#} $ or
$J(\Gamma ) = J(\Lambda )$ and $\Gamma $ is hereditary.
\end{lemma}

\bew
Since $ \Lambda /J(\Lambda )$ is a simple 
$\Lambda $-$\Lambda $-bimodule, also its dual $\Gamma / \Lambda $ is simple.
Now 
 $\Gamma \supseteq J(\Gamma ) + \Lambda \supseteq \Lambda $ and
$J(\Lambda ) = J(\Gamma ) \cap \Lambda $.
Therefore either $\Gamma = J(\Gamma ) + \Lambda \neq \Lambda $ and
$\Gamma /J(\Gamma ) \cong \Lambda /J(\Lambda )$ or
$J(\Gamma ) \subseteq \Lambda $ whence $J(\Gamma ) = J(\Lambda )$.
 In the latter case $\Gamma = \Id(J(\Gamma )) $  is hereditary.
In the first case  $\Lambda $ is not hereditary and
$J(\Gamma )$ is the unique maximal $2$-sided $\Gamma $-ideal 
in $\Gamma $.
Since $J(\Lambda ) \subset J(\Gamma) \subset J(\Lambda )^{\#} = \Gamma $
(by Theorem \ref{jistdual}),
one also has
$J(\Lambda ) \subset J(\Gamma)^{\# } \subset  \Gamma $.
Here all inclusions are proper.
So  $J(\Gamma )^{\# }$ is also a maximal $2$-sided
$\Gamma $-ideal in $\Gamma $
and hence $J(\Gamma )^{\# } = J(\Gamma )$.
\eb

Following the lines of \cite[11.4]{Jac} one gets the following
remark:

\begin{rem}
With the assumptions of Lemma \ref{selbstduallok}
assume that $J(\Lambda ) = J(\Gamma )$.
Then $\Gamma $ is hereditary and 
$\Gamma / J(\Gamma )$ has two (isomorphic) composition factors as
a $\Lambda - \Lambda $-bimodule, namely the submodule
$\Lambda / J(\Lambda) $ and its dual, the factor module
$\Gamma / \Lambda = (\Lambda / J(\Lambda )) ^{\#} $.
\end{rem}

\begin{notation}\label{notat}
We fix the following notation:
\\
$\Lambda $ denotes an indecomposable non-hereditary symmetric $R$-order in ${\cal A}$,
\\
 $\Gamma := \Id (J(\Lambda )) = J(\Lambda )^{\#} $ the idealizer of
the radical of $\Lambda $,
\\
$\epsilon _1,\ldots , \epsilon _s$ are the central primitive idempotents of
${\cal A}$, and 
\\
$e_1,\ldots , e_h \in \Lambda $ are orthogonal lifts of the central primitive 
idempotents of $\Lambda /J(\Lambda )$.
According to Lemma \ref{selbstduallok} we order the $e_i$ such that
$e_i\Lambda e_i/J(e_i \Lambda e_i) \cong e_i\Gamma e_i /J(e_i \Gamma e_i)$
for $1\leq i \leq t \leq h$ and 
$J(e_i\Gamma e_i ) = J(e_i \Lambda e_i)$ for $t < i \leq h $ and put
$$e:=\sum_{i=1}^t e_i \mbox{ and }
f:=\sum_{i=t+1}^h e_i.$$
\end{notation}

From Lemma \ref{selbstduallok} one now gets 

\begin{folg}{\label{efmoegl}}
For $1\leq i \leq h$ the idempotent
$e_i+J(\Gamma )$ is a central idempotent of $\Gamma / J(\Gamma )$.
If $1\leq i\leq t$, then $e_i + J(\Gamma )$ is a central primitive idempotent. 
If $t<i\leq h$ then $e_i \Gamma e_i $ is hereditary.
\end{folg}

\section{The second step of the radical idealizer chain}

We keep the Notation \ref{notat}. 
Moreover let $\Delta := \Id (J(\Gamma ))$ and 
let $\Lambda _N $ be the head order of $\Lambda $.

\begin{satz}{\label{ef}}
$$\Delta  = (J(e\Gamma e)^2+ e\Gamma f \Gamma e)^{\# } 
\oplus f\Gamma f \oplus f \Gamma e \oplus
e \Gamma f
.$$
\end{satz}

\bew
$\Delta = ((J(\Gamma )^{\# } J(\Gamma ) + J(\Gamma ) J(\Gamma )^{\#})^{\#}$
by Proposition \ref{ddd}.
Now
$$J(\Gamma )^{\#}  = (e J(\Gamma ) e)^{\#}  \oplus f\Gamma f \oplus
e\Gamma f \oplus f \Gamma e .$$
Lemma \ref{selbstduallok} says $(eJ(\Gamma )e)^{\#} = e J(\Gamma ) e$
and therefore
$$J(\Gamma )^{\#} J(\Gamma ) = 
(e J(\Gamma ) e  \oplus f\Gamma f \oplus e \Gamma f \oplus f \Gamma e) 
J(\Gamma ) =  $$
$$ ((eJ(\Gamma ) e)^2 + e\Gamma f \Gamma e) \oplus ( f J(\Gamma ) f +
f \Gamma e \Gamma f) \oplus (f\Gamma e) \oplus (eJ(\Gamma ) e \Gamma f
+ e\Gamma f J(\Gamma )f ) .$$
Since $$f\Gamma e \Gamma f = 
f \Lambda e \Lambda f \subseteq J(f\Lambda f) = J(f\Gamma f) $$ one gets
$$\Delta ^{\#} = 
((eJ(\Gamma ) e)^2 + e\Lambda f \Lambda e) \oplus  f J(\Lambda ) f 
\oplus f\Lambda  e \oplus e \Lambda f$$ and therefore 
the theorem follows.
\eb

\begin{prop}
For the head order $\Lambda _N$ one finds
$$ f\Lambda _N f = f \Gamma f ,\ e\Lambda _N f = e \Lambda f ,\ 
f\Lambda _N e = f\Lambda e  $$
and $e +J(\Lambda _N)$ and $f + J(\Lambda _N)$ are central idempotents of 
 $\Lambda _N / J(\Lambda _N)$.
\end{prop}

\bew
Let $\Lambda _0 := \Lambda $ and
$\Lambda _i := \Id(J(\Lambda _{i-1}))$
 ($1\leq i \leq N$).

Using induction we show that 
$ f\Lambda _i f = f \Gamma f ,\ e\Lambda _i f = e \Lambda f ,\ 
f\Lambda _i e = f\Lambda e $ and that
$f + J(\Lambda _i)$ (hence also $e+J(\Lambda _i)$) 
lies in the center of  $\Lambda _i/J(\Lambda _i)$
for all $1\leq i \leq N$.

For $i=1$ this is trivial.
Now let 
 $i\geq 1$  and  assume that the statement is true for $i$.

Then $\Lambda _{i+1} = 
((J(\Lambda _i)^{\# } J(\Lambda _i)+ J(\Lambda _i)J(\Lambda _i)^{\#})^{\#}$.
By assumption $$J(\Lambda _i ) = J(e\Lambda _i e) \oplus J(f\Lambda f) \oplus
e \Lambda f \oplus f\Lambda e $$ 
and therefore
$$J(\Lambda _i ) ^{\#}  = J(e\Lambda _i e) ^{\#}  \oplus f \Gamma f \oplus
e \Lambda f \oplus f\Lambda e .$$

Since  $\Lambda _i $ radically covers  $\Lambda _{i-1} $, it holds that
$J(\Lambda _{i-1} ) \subseteq J(\Lambda _i)$.
In particular $J(\Lambda _0) \subseteq J(\Lambda _i )$.
Therefore $J(\Lambda _i) ^{\#} \subseteq J(\Lambda _0)^{\#} = \Gamma $.
One calculates
$J(\Lambda _i) ^{\#}  J(\Lambda _i) + J(\Lambda _i) J(\Lambda _i) ^{\#} $
$$=(J(e\Lambda _i e)^{\#}J(e\Lambda _ie)+J(e\Lambda _ie)J(e\Lambda _i e)^{\#} 
 + e \Lambda f \Lambda e ) \oplus
J(f\Lambda f) \oplus e\Lambda f  \oplus f\Lambda e . $$
After dualizing, one gets the desired form of
$\Lambda _{i+1}$. 

Let $\lambda \in \Lambda _{i+1}$. Then
$$f\lambda - \lambda f = f \lambda e - e\lambda f \in
e\Lambda _{i+1} f \oplus f\Lambda _{i+1} e =
e\Lambda  f \oplus f\Lambda e \subseteq J(\Lambda )\subseteq J(\Lambda_{i+1}).$$
Therefore $f+J(\Lambda _{i+1}) \in Z(\Lambda _{i+1} / J(\Lambda _{i+1}))$
 and hence also  $e+J(\Lambda _{i+1})=(1-f)+J(\Lambda _{i+1})$
is central.
\eb

\begin{folg}
The conductor of $\Lambda _N$ in $\Lambda $ is
$$\Lambda _N^{\#} = (e\Lambda _N e)^{\#} \oplus f J(\Lambda ) f
\oplus e\Lambda f \oplus f \Lambda e .$$
\end{folg}

\begin{satz}{\label{f=01}}
Let $\Lambda $, $\Gamma =\Id(J(\Lambda ))$, $\Delta =\Id(J(\Gamma )) $
 and $f$ be as in Notation \ref{notat}.
\\
Then
$f=0$ or $f=1$.
\\
If $f=1$, then $J(\Gamma ) = J(\Lambda )$ and 
$\Gamma = \Delta $ is hereditary.
\\
If $f=0$, then $\Gamma = J(\Gamma ) + \Lambda $, 
$\Gamma / J(\Gamma ) \cong \Lambda /J(\Lambda ) $,
$J(\Gamma ) = J(\Gamma )^{\#} $  and 
$\Delta = (J(\Gamma )^2)^{\#}.$
\end{satz}

\bew
With Notation  \ref{condprec}
 the conductor formula \ref{condform} gives
$$\Lambda _N^{\#} = (e\Lambda _N e)^{\#} \oplus f J(\Lambda ) f
\oplus e\Lambda f \oplus f \Lambda e  =
\oplus_{i=1}^s  \wp _i^{m_i(n_i-d_i-1)} 
( eJ(\Lambda _N) e \oplus f J(\Lambda ) f
\oplus e\Lambda f \oplus f \Lambda e ) .$$
So if $\epsilon _i f \neq 0$ for some $1\leq i \leq s$, then
$n_i-d_i=1$ and $J(\Lambda _N)\epsilon _i \subseteq \Lambda $.
Since  $J(\Lambda _N)\epsilon _i $ is a  pro-nilpotent $\Lambda $-ideal in
$\Lambda$ containing
 $J(\Lambda )\epsilon _i$ one gets
 $J(\Lambda ) \epsilon _i = J(\Lambda _N) \epsilon _i \subseteq J(\Lambda )$.
But then $\epsilon _i \in \Gamma $ and $\Gamma \epsilon _i =
\Lambda _N \epsilon _i$ is hereditary.

We claim that $f\epsilon _i = \epsilon _i$.
To see this let $1\leq j\leq t$ with $e_j \epsilon _i \neq 0$.
Since $e_j + J(\Gamma )$ is a central primitive idempotent of
$\Gamma / J(\Gamma )$, one even has $e_j \epsilon _i = e_j$.
Now $\Gamma \epsilon _i$ is hereditary, so 
$e_j\Gamma e_j \cong \Omega _i^{x\times x}$ for some $x\in \N$.
Let $P_i$ denote the maximal ideal  in $\Omega _i$. 
Since $j\leq t$
Lemma  \ref{selbstduallok} says that 
$P _i ^{x\times x} \cong J(e_j \Gamma e_j)$ is symmetric with respect to
the restriction of the form $Tr_z$  above.
But $n_i-d_i-1=0$, yields
together with \cite[Theorem 14.9]{Rei} that
$J(e_j \Gamma e_j)^{\#} = \wp _i^{-1} e_j\Gamma e_j$ which is a contradiction.
 Therefore $e \epsilon _i = 0$ and hence
$f \epsilon _i = \epsilon _i$.

So for all central primitive idempotents $\epsilon _i$ of ${\cal A}$
either $f\epsilon _i =0 $  or $f\epsilon _i = \epsilon _i$.
Therefore $\Lambda = e\Lambda e \oplus
f\Lambda f$. Since  $\Lambda $ is assumed to be indecomposable one has 
$f=0$ or $f=1$.
In the latter case $J(\Gamma ) = J(\Lambda ) = J(\Lambda _N)$,
hence $\Gamma $ is hereditary.
If $f=0$, then $\Delta = (J(\Gamma )^2)^{\# }$ from Theorem \ref{ef}.
The fact that $J(\Gamma ) $  is self-dual follows with Lemma \ref{selbstduallok}
which also implies that $\Gamma  = J(\Gamma ) + \Lambda $ and
hence $\Lambda /J(\Lambda ) \cong \Gamma /J(\Gamma )$.
\eb

Summarizing let $\Lambda $ be an indecomposable, non-hereditary,
symmetric $R$-order in ${\cal A}$, $\Gamma = \Id(J(\Lambda ))$ and
$\Delta = \Id(J(\Gamma ))$. 
Let $f$ be as in Notation \ref{notat}.
Then 

\begin{picture}(400,140)(0,0)

\put(50,125){\circle*{4}}
\put(50,100){\circle*{4}}
\put(50,75){\circle*{4}}
\drawline(50,125)(50,75)
\drawline(300,125)(325,100)(300,75)
\drawline(300,125)(275,100)(300,75)
\put(300,125){\circle*{4}}
\put(325,100){\circle*{4}}
\put(275,100){\circle*{4}}
\put(300,75){\circle*{4}}
\put(48,125){\makebox(0,0)[r]{$\Gamma $}}
\put(48,100){\makebox(0,0)[r]{$\Lambda $}}
\put(48,75){\makebox(0,0)[r]{$J(\Lambda ) = J(\Gamma )$}}

\put(50,50){\makebox(0,0)[c]{$f=1$}}
\put(50,25){\makebox(0,0)[c]{$\Gamma $ hereditary}}

\put(302,125){\makebox(0,0)[l]{$\Gamma $}}
\put(327,100){\makebox(0,0)[l]{$\Lambda $}}
\put(273,100){\makebox(0,0)[r]{$J(\Gamma ) = J(\Gamma )^{\#}$}}
\put(302,75){\makebox(0,0)[l]{$J(\Lambda ) $}}

\put(300,50){\makebox(0,0)[c]{$f=0$}}
\put(300,30){\makebox(0,0)[c]{$\Gamma /J(\Gamma ) \cong \Lambda/J(\Lambda )$}}
\put(300,15){\makebox(0,0)[c]{$\Delta = \Id (J(\Gamma ))  = (J(\Gamma )^2)^{\#} $}}
\put(150,75){\makebox(0,0)[c]{or}}

\end{picture}

Returning to the proof of Theorem \ref{ef} with the two possibilities
$f=1$ or $f=0$ we find that either $\Gamma = \Delta = 
O_r(J(\Lambda ))$ is hereditary or 
$J(\Gamma ) = J(\Gamma )^{\#}$ whence $ O_r(J(\Gamma )) =
(J(\Gamma ) ^{\#} J(\Gamma ) )^{\#}  = \Delta $.
Therefore the first two steps in the right- and left-radical idealiser
chain of a symmetric order $\Lambda $ coincide.

\begin{kor}
Let $\Lambda $ be a symmetric order, 
$\Gamma = \Id (J(\Lambda )) $ (which equals $ O_r (J(\Lambda )) = O_l(J(\Lambda ))$ by Corollary \ref{l=r}).
Then  
$$ \Id (J(\Gamma )) = O_l (J(\Gamma )) = O_r (J(\Gamma )) .$$
\end{kor}

\section{Symmetric orders with radical idealizer length 1 or 2.}

The case $f=1$ in Theorem \ref{f=01} can be dealt with the arguments 
in \cite[Satz 11.4]{Jac}.
With a modest additional assumption one gets $J(\Gamma ) = J(\Lambda )$
and in particular if $\Lambda $ is a block of a group ring with
$l_{rad}(\Lambda ) = 1$ one can associate a Brauer tree to $\Lambda $.

\begin{prop}{\label{j=j}}
Assume that $\Gamma = \Id (J(\Lambda ))$ is hereditary.
Assume further that 
 $Z_i := \epsilon _i Z(\Lambda ) \subseteq \epsilon _i {\cal A} $ 
is a maximal order
for all central  primitive idempotents $\epsilon _1,\ldots , \epsilon _s$.
Then $$J(\Gamma ) = J(\Lambda ).$$
\end{prop}

\bew
With the Notation \ref{condprec} and 
 Jacobinski's conductor formula \ref{condform}
$$J(\Lambda ) = \Gamma ^{\#} = \bigoplus _{i=1}^s \epsilon _i 
\wp  _i ^{m_i(n _i - \delta _i - 1)} J(\Gamma ) .$$
By Remark \ref{piinJ} and since $\epsilon _i J(\Lambda ) \subset J(\Lambda )$
$$J(Z_i) J(\Gamma ) \subsetneq J(Z_i) \Gamma \subseteq J(\Lambda )$$
for all $i$. 
This implies that $n _i - \delta _i = 1$ and that $J(\Lambda ) = J(\Gamma )$.
\eb


\begin{prop}\label{tree}
Assume that $J(\Gamma ) = J(\Lambda ).$
If moreover 
 the decomposition map from the Grothendieck groups of 
simple modules $G_0({\cal A}) \to G_0(k\otimes \Lambda )$ is surjective
or $k$ is a splitting field for $k\otimes \Gamma $,
then for each idempotent $e_i \in \Lambda $, there are exactly 
two central primitive idempotents $\epsilon_{i_1} $ and $\epsilon _{i_2}$
in ${\cal A}$ with $e_i \epsilon _{i_j} \neq 0 $ ($j=1,2$).
\end{prop}

\bew
We make precise the embedding $\Lambda / J(\Lambda ) \hookrightarrow \Gamma / J(\Gamma )$ following the lines of the proof of \cite[Satz 11.4]{Jac}:
Let $e_i(\Lambda / J(\Lambda )) =:  S_i $ be the simple algebra summand of 
$\Lambda /J(\Lambda )$ that corresponds to $e_i$ ($1\leq i \leq h $).
Then $$\Lambda / J(\Lambda ) \cong S_1\oplus \ldots \oplus S_h \mbox{ and }
 \Gamma / \Lambda = (\Lambda  / J(\Lambda ))^{\#} \cong 
S_1^* \oplus \ldots \oplus S_h ^* \cong S_1 \oplus \ldots \oplus S_h $$
as $\Lambda - \Lambda $-bimodules.
Therefore for every simple $k$-algebra summand $S_i$ of $\Lambda /J(\Lambda )$
there are either two algebra-summands $T_{i_1'}$ and $T_{i_2'}$ of 
$\Gamma /J(\Gamma )$ such that $S_i$ is diagonally embedded into
$T_{i_1'} \oplus T_{i_2'}$ 
or there is a unique summand $T_i \cong l_2^{n\times n}$ such that
$S_i \cong l_1^{n\times n} \subset T_i$ for extension fields 
$l_2,l_1$ of $k$ with $[l_2:l_1] = 2$.
The latter is impossible, if $k$ is a splitting field for $k \otimes \Gamma $.
Similarly, if the decomposition map of $\Lambda $ is surjective, the last case 
cannot happen, since otherwise the simple $S_i$-module occurs with
even multiplicity in every simple $\Gamma $-module and hence in the
reduction of every $\Gamma $-lattice modulo $\pi $.
Therefore $$S_i \hookrightarrow T_{i_1} \oplus T_{i_2} \subset \Gamma / J(\Gamma ).$$
Since $\Gamma $ is hereditary, it contains the central primitive idempotents
$\epsilon _i$ of ${\cal A}$.
Therefore $\Gamma / J(\Gamma ) = \oplus _{i=1}^s \epsilon _i \Gamma / J(\epsilon _i \Gamma ) $ and hence each simple summand of 
$\Gamma / J(\Gamma )$ is a summand of some 
$ \epsilon _i \Gamma / J(\epsilon _i \Gamma ) $.
In particular for $j=1,2$ the summand 
$T_{i_j'}$  defines a unique central primitive
idempotent $\epsilon _{i_j}$ with 
$f_{i_j'} \epsilon _{i_j}   \neq 0$ for any lift $f_{i_j'} \in \Gamma $ of
the central primitive idempotent of $\Gamma / J(\Gamma )$ that belongs to $T_{i_j'}$.
Then $\epsilon _{i_j} $ ($j=1,2$) are the only central primitive idempotents 
in ${\cal A}$ with $e_i \epsilon _{i_j} \neq 0$.
\eb

If the decomposition map of $\Lambda $ is surjective 
(which is always satisfied when $\Lambda $ is a block of a group ring)
or $k$ is a splitting field for $k\otimes \Gamma $ and 
the other  assumptions of Proposition \ref{tree} hold, we can define a graph
${\cal G}(\Lambda )$  whose vertices correspond to 
$\epsilon _1,\ldots ,\epsilon _s$ and whose edges correspond to 
$e_1,\ldots , e_h$.
Two vertices $\epsilon _i$ and $\epsilon _j$ are connected 
by the edge $e_l$, if $e_l \epsilon_i \neq 0 $ and $e_l \epsilon _j \neq 0$.

As in \cite[Korollar 11.6]{Jac} one shows:
\begin{folg}
If the decomposition map of $\Lambda $ is surjective then
${\cal G}(\Lambda )$ is a tree.
\end{folg}

We end this section with a short remark on the length 2 case.

%
%

\begin{rem}
With Notation \ref{notat}
assume that $f=0$ and $\Delta :=\Id (J(\Gamma ))$ is hereditary.
\bi
\item[(i)]
$J(\Lambda )^2 \epsilon _i \in \Lambda $
for all $1\leq i \leq s $. In particular 
$\pi ^2 \epsilon _i \in \Lambda $ which means that the defect of
$\Lambda $ is  $\leq 2$.
\item[(ii)]
$\epsilon _i e_l\Lambda e_j \subseteq e_l \Lambda e_j $ for all 
$1\leq i \leq s $, $1\leq l\neq j \leq h$.
\item[(iii)]
If $s>1$ then
$\epsilon _i J(\Lambda ) \not\subseteq J(\Lambda )$ for all $1\leq i \leq s$.
\ei
\end{rem}

\bew
(i)
Let $1\leq i \leq s$. Since $\epsilon _i \in \Delta $ one has
$\epsilon _i J(\Gamma ) \subseteq  J(\Gamma )$.
Now $$J(\Lambda ) ^2  \subseteq J(\Gamma )^2 = \Delta ^{\# } =\bigoplus _{i=1}^s J(\Gamma )^2 \epsilon _i \subseteq
\Lambda $$ implies $J(\Lambda ) ^2 \epsilon _i \in \Lambda $.
Since $\pi = \pi 1 \in \pi \Lambda \subset J(\Lambda )$ one
has $\pi^2 \epsilon _i \in \Lambda $ for all $i$.
\\
(ii)
Let $l\neq j \in \{ 1,\ldots , h\}$. Then
$e_l\Lambda e_j \subseteq J(\Gamma )$.
Since $\epsilon _i J(\Gamma ) \subseteq J(\Gamma )$, the claim follows.
\\
(iii)
Assume that there is $1\leq i \leq s$ such that $\epsilon _i J(\Lambda ) 
\subseteq J(\Lambda )$.
Then $\epsilon _i \in \Gamma $.
Since $\epsilon _i \neq 1 $ and $\Lambda $  is indecomposable there is
$1\leq j\leq h$ such that $0\neq \epsilon _i e_j \neq e_j$.
But then $e_j+J(\Gamma ) \in Z(\Gamma / J(\Gamma ))$ is not primitive
contradicting Corollary \ref{efmoegl}.
\eb


\section{$p$-groups}

Let $G\neq \{ 1 \}$ be a  $p$-group, $R=\Z _p$,  and $\Lambda := RG$.
Then $\Lambda $ is a symmetric $R$-order with respect to the associative 
bilinear form 
$$\phi (x,y):=  \frac{1}{|G|} \tr _{reg} (xy)  =(xy) _1 \mbox{ if } xy = \sum _{g \in G} (xy)_g g $$
where $\tr _{reg}$ is the regular trace of $\Q_p G$.
Moreover
$$ J(\Lambda ) = \langle p \Lambda , g-h \mid g,h\in G \rangle _R  $$
and 
$$\Gamma := \Id (J(\Lambda )) = J(\Lambda ) ^{\# } 
= \langle \Lambda , \frac{1}{p} \sum _{g\in G} g \rangle _R .$$
because $\frac{1}{p} \sum _{g\in G} g $ idealizes $J(\Lambda )$ 
and $\Gamma $ is an over-order of $\Lambda $ of index $p = |\Lambda /J(\Lambda )|$.
If $|G|=p$ then $\Gamma $ is hereditary by \cite[Section 11]{Jac}.
Therefore we assume that $|G| \geq p^2$.
Then the radical of $\Gamma $ is
$$J(\Gamma ) = \langle J(\Lambda ) , \frac{1}{p} \sum _{g\in G} g \rangle _R 
= J(\Gamma )^{\# } $$
and is contained in $\Gamma $ of index $p$.

\begin{theorem}\label{id2p}
Let $\Delta := \Id (J(\Gamma )) $. Then
$$\Delta = (J(\Gamma )^2)^{\#} =
\langle \Lambda, \frac{1}{p^2} \sum _{g\in G} g, \frac{1}{p} \sum _{g\in G} \varphi (g) g \mid
\varphi \in \Hom ( G , R/pR) \rangle _R  .$$
\end{theorem}

\bew
Clearly $\lambda :=\frac{1}{p^2} \sum _{g\in G} g \in \Id (J(\Gamma )) $.
Let $y:=\sum _{g\in G} y_g g \in \Id (J(\Gamma ))$.
Then $$p \sum _{g\in G} y_g g = a \frac{1}{p} \sum _{g\in G} g + x $$
with $x\in J(\Lambda ) \subset \Lambda $ and $a\in R$.
Hence $p y_g \equiv \frac{a}{p} \pmod{R} $.
Replacing $y$ by $y-a\lambda $ we may assume that 
$y_g = \frac{a_g}{p}$ with $a_g\in R$ for all $g\in G$.
Adding a suitable multiple of $p\lambda  $ to $y$ we can also assume that 
$a_1 = 0$.
Now $\Id (J(\Gamma )) = (J(\Gamma )^2)^{\# }$ and one calculates
$$J(\Gamma )^2 = \langle p^2 \Lambda , p(g-h),(g_1-h_1)(g_2-h_2) , \sum_{g\in G} g \mid g,h,g_1,h_1,g_2,h_2 \in G \rangle _R .$$
In particular the coefficient of $1$ of 
$y (g^{-1}-1)(1-h^{-1}) $ which is 
$\frac{1}{p} (a_g+a_h-a_{hg})$ lies in $R$.
Hence $$\varphi : g \mapsto a_g + pR \in R/pR $$ 
is a group homomorphism from $G$ to $R/pR$, from which the inclusion 
$\subseteq $ follows.
It remains to show that the elements 
$y := \frac{1}{p} \sum _{g\in G } \varphi (g) g $ with
$\varphi \in \Hom (G, R/pR)$ are in the dual of $J(\Gamma )^2$.
Clearly 
$$\phi (y,p^2\Lambda ) \subset R \mbox{ and } \phi (y,p(g-h)) \subset R \fa g,h \in G .$$
For $(g_1-h_1)(g_2-h_2)$ with $g_1,g_2,h_1,h_2\in G $ one gets 
$$\phi (y,
(g_1-h_1)(g_2-h_2) )= \frac{1}{p} (\varphi (g_2^{-1}g_1^{-1} ) 
-\varphi (h_2^{-1} g_1^{-1} ) - \varphi (g_2^{-1} h_1^{-1}) + \varphi (h_2^{-1} h_1^{-1} )) \in R $$
since $\varphi $ is a homomorphism.
The last generator is $p^2\lambda $ for which one finds
$$\phi (y,\sum _{g\in G} g ) = \frac{1}{p} \sum _{g\in G} \varphi (g) \in R.$$
\eb

\begin{folg}\label{index}
Let $\Lambda := \Z_pG$ for some $p$-group $G$ of order $|G| \geq p^2$.
Let $\Gamma := \Id(J(\Lambda )) $ and $\Delta := \Id (J(\Gamma ))$.
Then
$|\Gamma / \Lambda | = p$ and
$|\Delta / \Lambda | = p^2| G/(G' G^p) |$.
\end{folg}

\begin{folg}
Let $G$ be a $p$-group and $\Lambda := \Z_p G$.
\begin{itemize}
\item[1)] $\Gamma = \Id (J(\Lambda ) ) $ is hereditary if and only if 
$G \cong C_p $.
\item[2)]
Assume that $|G| \geq p^2$. Then
$\Delta = \Id (J(\Gamma )) $ is hereditary if and only if $|G| =4$.
\end{itemize}
\end{folg}

\bew
1) Since $\Gamma \subset \frac{1}{p} \Lambda $ contains all central 
primitive idempotents of $\Q _p G $, one gets $|G| = p$.
That $\Id(J(\Z_pC_p)) $ is hereditary follows from 
 \cite[Abschnitt 11]{Jac}.
\\
2) Since $\Delta \subset \frac{1}{p^2} \Lambda $ contains all central 
primitive idempotents of $\Q _p G $, one gets $|G| = p^2$.
Hence $G$ is Abelian, $G\cong C_{p^2}$ or $G\cong C_p\times C_p $ and
by Corollary \ref{index}, $[\Delta : \Lambda ] = p^3$ respectively $p^4$.
If $\Delta $ is hereditary, then $\Delta $ is the maximal order in
$\Q_p G$, 
$$\Delta \cong \Z_p \oplus \Z_p [\zeta _p ] \oplus \Z_p [\zeta _{p^2}] 
\mbox{ respectively }
\Delta \cong \Z_p \oplus \Z_p [\zeta _p ] ^{p+1} .$$ 
The form $\phi $ above is $\frac{1}{|G|} \tr _{reg} $ hence
the discriminant of $\Delta $ with respect to 
$\phi $ is 
$$p^{-2p^2}  \cdot 1 \cdot p^{p-2} \cdot p^{p(2p-3)}  = p^{-2p-2}
\mbox{ respectively }
p^{-2p^2} \cdot 1 \cdot (p^{p-2})^{p+1} = p^{-p^2-p-2} $$
(see \cite[Prop. 2.1]{Was}).
Since $\Delta $ contains a symmetric order of index $p^3$ respectively
$p^4$, one gets $-2p-2 = -6 $ hence $p=2$ respectively
$-p^2-p-2=8$ which also implies $p=2$.
The same argument shows that for $G= C_2\times C_2$ or $G=C_4$, 
the order $\Delta $ has the same discriminant as the maximal order
and hence is hereditary (i.e. equal to the maximal order).
\eb

Note that this corollary also 
follows from Theorem \ref{abelian} below.

For Abelian groups $G$, the radical idealizer length of 
$\Z_p G$ can be calculated from the exponent and the order
 of the Sylow $p$-subgroup of $G$:

\begin{theorem}\label{abelian}
Let $G$ be an Abelian group with Sylow $p$-subgroup 
of order $p^n$ and of exponent $p^a>1$.
Then 
$$l_{rad} (\Z_pG) = p^{a-1} + (p^{a} - p ^{a-1}) (n-a) .$$
\end{theorem}

\proof
The theorem follows with Corollary \ref{folgcommut} by calculating the conductor of 
the maximal order $\Gamma$
in $\Q_pG$:
$$\Gamma = 
\bigoplus _{i=1}^{s} \Z_pG \epsilon _i =
\bigoplus _{i=1}^{s} R_i [\zeta _{p^{a_i}}] $$
where $R_i$ is an unramified extension of $\Z_p$ and $\{a_1,\ldots , a_s\} = 
\{ 0,1,\ldots , a \} $.
If $\ ^*$ denotes the different, i.e. the dual with respect to the
usual trace bilinear form, then by 
\cite[Prop. 2.1]{Was} 
$ R_i [\zeta _{p^{a_i}}]  ^{*} =
 R_i [\zeta _{p^{a_i}}]  (1-\zeta _{p^{a_i}})^{ -p^{a_i-1}(a_ip-a_i-1) }$
and hence the conductor of $\Gamma $ in $\Z_pG $ is 
$$\Gamma ^{\#} = 
\bigoplus _{i=1}^{s} R_i [\zeta _{p^{a_i}}]  ^{*} p^n =
\bigoplus _{i=1}^{s} R_i [\zeta _{p^{a_i}}] 
 (1-\zeta _p^{a_i})^{ ((n-a_i) (p-1) +1)p^{a_i-1}} .$$
By Corollary \ref{folgcommut}, the length of the radical idealizer chain is 
$$\max _{i=0,\ldots , a} 
 ((n-i) (p^i-p^{i-1}) + p^{i-1} )
 = (n-a) (p^a-p^{a-1}) + p^{a-1} 
.$$
\eb

\renewcommand{\arraystretch}{1}
\renewcommand{\baselinestretch}{1}
\large

\end{document}